\documentclass[10pt,twoside]{article}
\usepackage{mathrsfs}
\usepackage{amsmath}
\usepackage{amssymb}
\usepackage{fancyhdr}
\usepackage{latexsym}
\usepackage{bbding}
\usepackage{mathrsfs}
\usepackage{wasysym}
\usepackage{cite}

\usepackage{multicol,graphics}

\newtheorem{theorem}{Theorem}[section]
\newtheorem{lemma}[theorem]{Lemma}

\newtheorem{definition}[theorem]{Definition}

\numberwithin{equation}{section}

\allowdisplaybreaks

\begin{document}
\title{{Global existence of large solutions for the parabolic-elliptic Keller--Segel system in Besov type spaces}
\thanks{This work is partially supported by the National Natural Science Foundation of China (nos. 11961030 and 12361034) and
 the Natural Science Foundation of Shaanxi Province (no. 2022JM-034).} }

\author{
 {\small   Jihong Zhao\thanks{\text{E-mail address}:
jihzhao@163.com.}}
\\
 {\small    School of Mathematics and Information Science, Baoji University of Arts and Sciences,}\\
{\small   Baoji, Shaanxi 721013, China} }
\date{}
\maketitle

\begin{abstract}
In this paper, we investigate global existence of large solutions for the parabolic-elliptic Keller--Segel system in the homogeneous Besov type spaces. A class of initial data was presented, generating a global smooth solution although the $\dot{B}^{-2}_{\infty,\infty}$-norm of the initial data may be chosen arbitrarily large.
\end{abstract}
\smallbreak

\medskip


\section{Introduction}

\setcounter{section}{1}
\setcounter{theorem}{0}

In this paper, we consider the global existence of large solutions for the Cauchy problem of the parabolic-elliptic Keller--Segel system:
\begin{equation}\label{eq1.1}
\begin{cases}
  \partial_{t}u-\Delta u+\nabla\cdot(u\nabla v)=0,\\
    -\Delta v=u,\\
  u(x,0)=u_0(x),
\end{cases}
\end{equation}
where two unknown functions $u$ and $v$ representing the cell density and the concentration of the chemical attractant, respectively, and both are unknown functions of the space variable $x\in \mathbb{R}^{d}$ and of the time variable $t\in\mathbb{R}^{+}$.

The system \eqref{eq1.1} is a simplified version of the mathematical model of chemotaxis, which was first introduced by E.F. Keller and L.A. Segel \cite{KS70} in 1970. This basic model
was used to describe the collective movement of bacteria possibly leading to cell aggregation by chemotactic
effect. It is also connected with astrophysical models of
gravitational self-interaction of massive particles in a cloud or a nebula,  see \cite{BHN94B,BCGK04}.  We refer to \cite{B98,  FP11,GZ98, HV96, JL92, KS08, L13, LW23,OS10, RZ95,W10, W13} for more details concerning about the existence and nonexistence (blow-up) of global in time solutions to the different Keller--Segel
models studied since the 1970s.

Our aim here is to show the global existence of large solutions of the system \eqref{eq1.1} in critical Besov type spaces. Karch \cite{K99} established global well-posedness with small initial data in critical Besov space $\dot{B}^{-2+\frac{d}{p}}_{p, \infty}(\mathbb{R}^d)$ with $\frac{d}{2}<p<d$. Subsequently,  Iwabuchi \cite{I11} proved local well-posedness and global well-posedness for small initial data in critical Besov spaces $\dot{B}^{-2+\frac{d}{p}}_{p,q}(\mathbb{R}^{d})$ and Fourier--Herz spaces $\dot{\mathcal{B}}^{-2}_{2}(\mathbb{R}^{d})$ with $1\leq p<\infty$ and $1\leq q\leq\infty$, moreover, the ill-posedness is also established in critical Besov spaces $\dot{B}^{-2}_{\infty,q}(\mathbb{R}^{d})$ with $2<q\leq\infty$. Finally, Iwabuchi--Nakamura \cite{IN13} showed   global existence of solutions for small initial data in $\dot{F}^{-2}_{\infty,2}(\mathbb{R}^{d})=BMO^{-2}$ through the Triebel--Lizorkin spaces, which combining the well-posed results in \cite{I11},  we know that, the simplified Keller--Segel system \eqref{eq1.1} is well-posed in $\dot{B}^{-2+\frac{d}{p}}_{p,\infty}(\mathbb{R}^{d})$ and $\dot{F}^{-2}_{\infty,2}(\mathbb{R}^{d})$ with $\max\{1, \frac{d}{2}\}<p<\infty$, but ill-posed in $\dot{B}^{-2}_{\infty,q}(\mathbb{R}^{d})$ with $2<q\leq\infty$. For more well-posed and ill-posed results to more general chemotaxis models, we refer the readers to see \cite{DL17, LYZ23, NY20, NY 22, XF22}.

Before we state the main results, let us first introduce the definition of the homogeneous Besov spaces.

\begin{definition}\label{de1.1}
Let $\varphi\in\mathcal{S}(\mathbb{R}^{d})$ be a positive radial function such that  $\varphi$ is supported in the shell $\mathcal{C}=\{\xi\in\mathbb{R}^{d},\ \frac{3}{4}\leq
|\xi|\leq\frac{8}{3}\}$, and
$
   \sum_{j\in\mathbb{Z}}\varphi(2^{-j}\xi)=1 \ \ \text{for}\ \  \xi\in\mathbb{R}^{d}\backslash\{0\}.
$
Let $h=\mathcal{F}^{-1}\varphi$. Then  for any $f\in\mathcal{S}'(\mathbb{R}^{d})$, we define the homogeneous dyadic block $\Delta_{j}$ as
\begin{align*}
 \Delta_{j}f(x): =2^{dj}\int_{\mathbb{R}^{d}}h(2^{j}y)f(x-y)dy,
\end{align*}
and for any $s\in \mathbb{R}$, $1\leq p,r\leq\infty$, we define the homogeneous Besov space $\dot{B}^{s}_{p,r}(\mathbb{R}^{d})$ by
\begin{equation*}
  \dot{B}^{s}_{p,r}(\mathbb{R}^{d}):=\Big\{f\in \mathcal{S}'(\mathbb{R}^{d})/\mathcal{P}:\ \
  \|f\|_{\dot{B}^{s}_{p,r}}<\infty\Big\},
\end{equation*}
\begin{equation*}
  \|f\|_{\dot{B}^{s}_{p,r}}:= \begin{cases} \left(\sum_{j\in\mathbb{Z}}2^{srj}\|\Delta_{j}f\|_{L^{p}}^{r}\right)^{\frac{1}{r}}
  \ \ &\text{for}\ \ 1\leq r<\infty,\\
  \sup_{j\in\mathbb{Z}}2^{sj}\|\Delta_{j}f\|_{L^{p}}\ \
  &\text{for}\ \
  r=\infty,
 \end{cases}
\end{equation*}
where $\mathcal{P}$ denotes the set of all polynomials.
\end{definition}

Next, we recall the definition of the so-called Chemin--Lerner mixed time-space spaces.

\begin{definition}\label{de1.2}  For $0<T\leq\infty$, $s\in \mathbb{R}$ and
$1\leq p, r, \rho\leq\infty$. We define the mixed time-space $\mathcal{L}^{\rho}(0,T; \dot{B}^{s}_{p,r}(\mathbb{R}^{d}))$
as the completion of $\mathcal{C}([0,T], \mathcal{S}(\mathbb{R}^{d}))$ by the norm
$$
  \|f\|_{\mathcal{L}^{\rho}_{T}(\dot{B}^{s}_{p,r})}:=\left(\sum_{j\in\mathbb{Z}}2^{srj}\left(\int_{0}^{T}
  \|\Delta_{j}f(\cdot,t)\|_{L^{p}}^{\rho}dt\right)^{\frac{r}{\rho}}\right)^{\frac{1}{r}}<\infty
$$
with the usual change if $\rho=\infty$
or $r=\infty$.
\end{definition}

Now we are ready to state our main results.  The first result is a new global existence result, under a nonlinear smallness assumption on the initial data.

\begin{theorem}\label{th1.1}  For any $u_{0}\in
\dot{B}^{-2+\frac{d}{p}}_{p,1}(\mathbb{R}^{d})$ with $1< p<\infty$, there exist two
positive constants $c_{0}$ and $C_{0}$ such that if
the initial data $u_{0}$ satisfies
\begin{align}\label{eq1.2}
  C_{0}\|e^{t\Delta}u_{0}\nabla(-\Delta)^{-1} e^{t\Delta}u_{0}\|_{L^{1}(\mathbb{R}^{+}; \dot{B}^{-1+\frac{d}{p}}_{p,1})}\exp\{C_{0}\|u_{0}\|_{\dot{B}^{-2+\frac{d}{p}}_{p,1}}
  \}\leq c_{0},
\end{align}
then the system \eqref{eq1.1} admits a unique global
solution $u$ associated with $u_{0}$, satisfying
\begin{equation}\label{eq1.3}
  u\in \mathcal{L}^{\infty}(\mathbb{R}^{+};
  \dot{B}^{-2+\frac{d}{p}}_{p,1}(\mathbb{R}^{d}))\cap  L^{1}(\mathbb{R}^{+};
  \dot{B}^{\frac{d}{p}}_{p,1}(\mathbb{R}^{d})).
\end{equation}
Moreover, there exists a positive constant $\kappa$ such that
\begin{equation}\label{eq1.4}
 \|u-e^{t\Delta}u_{0}\|_{\mathcal{L}^{\infty}(\mathbb{R}^{+}; \dot{B}^{-2+\frac{d}{p}}_{p,1})}
  +\kappa\|u-e^{t\Delta}u_{0}\|_{L^{1}(\mathbb{R}^{+}; \dot{B}^{\frac{d}{p}}_{p,1})}\leq
  \frac{c_{0}}{2}.
\end{equation}
\end{theorem}

The proof of Theorem \ref{th1.1} is given in Section 2 below. It consists in writing the solution $u$ (which exists for a short
time at least), as $u=e^{t\Delta}u_{0}+\bar{u}$ and in proving a global existence result for the perturbed Keller--Segel system satisfied by $\bar{u}$, under assumption \eqref{eq1.2}.  Moreover, the second result of this paper is to exhibit examples of applications of Theorem \ref{th1.1} which go beyond the
assumption of smallness of $\|u_{0}\|_{\dot{B}^{-2}_{\infty,\infty}}$-norm. For the sake of simplicity, we take $d=3$ to construct such an initial data.

\begin{theorem}\label{th1.2}  Let $\phi\in\mathcal{S}(\mathbb{R}^3)$ be a given function, and consider two real numbers $\varepsilon$, $\alpha\in(0,1)$. Define
 \begin{equation*}
 \varphi_{\varepsilon}(x):=\frac{(-\log \varepsilon)^{\frac{1}{5}}}{\varepsilon^{1-\alpha}}\cos(\frac{x_{3}}{\varepsilon})\phi(x_{1}, \frac{x_{2}}{\varepsilon^{\alpha}},x_{3}).
\end{equation*}
There exists a constant $C>0$ such that for $\varepsilon$ small enough, the smooth function $u_{0,\varepsilon}:=\nabla\cdot w_{0,\varepsilon}$ satisfies
\begin{equation}\label{eq1.5}
\|u_{0, \varepsilon}\|_{\dot{B}^{-2}_{\infty,\infty}}=\|\nabla\cdot w_{0, \varepsilon}\|_{\dot{B}^{-2}_{\infty,\infty}}\approx\|w_{0, \varepsilon}\|_{\dot{B}^{-1}_{\infty,\infty}}\geq  C^{-1}(-\log \varepsilon)^{\frac{1}{5}},
\end{equation}
where  $w_{0,\varepsilon}$ is defined as
$
w_{0,\varepsilon}:=(\partial_{2}\varphi_{\varepsilon}, \partial_{1}\varphi_{\varepsilon}, 0).
$
Moreover, we have
\begin{equation}\label{eq1.6}
  \|e^{t\Delta}u_{0, \varepsilon}\nabla (-\Delta)^{-1}e^{t\Delta}u_{0, \varepsilon}\|_{L^{1}(\mathbb{R}^{+}; \dot{B}^{-1+\frac{d}{p}}_{p,1})}\leq C\varepsilon^{\frac{\alpha}{3}}(-\log \varepsilon)^{\frac{2}{5}}.
\end{equation}
\end{theorem}

Notice that, for initial data $u_{0,\varepsilon}$ given by Theorem \ref{th1.2}, we can easily calculate that
\begin{align*}
  u_{0,\varepsilon}=\nabla\cdot w_{0,\varepsilon}=\frac{2}{\varepsilon}(-\log \varepsilon)^{\frac{1}{5}}\cos(\frac{x_{3}}{\varepsilon})\partial_{1}\partial_{2}\phi(x_{1} \frac{x_{2}}{\varepsilon^{\alpha}},x_{3}),
\end{align*}
and
\begin{align*}
\|u_{0,\varepsilon}\|_{L^{1}}\leq C \varepsilon^{\alpha}(-\log \varepsilon)^{\frac{1}{5}}\|\partial_{1}\partial_{2}\phi\|_{L^{1}}.
\end{align*}
Based on this fact, we know that, although the  $\dot{B}^{-2}_{\infty,\infty}$-norm of  $u_{0,\varepsilon}$ can be arbitrary large, but the $L^{1}$-norm of $u_{0,\varepsilon}$ tends to 0 as $\varepsilon$ tends to $0$. Hence, such initial data given by Theorem \ref{th1.2} can not lead to blow-up of solutions for the system \eqref{eq1.1}.

The proof of Theorem \ref{th1.2} will be given in Section 3. Throughout the paper, $C$ stands for a real positive constant which may be different in each occurrence.

\section{Proof of Theorem \ref{th1.1}}

Let us denote $u_{L}:=e^{t\Delta}u_{0}$, $\bar{u}:=u-u_{L}$. Then we can rewrite system \eqref{eq1.1} as
\begin{equation}\label{eq2.2}
\begin{cases}
\partial_{t}\bar{u}-\Delta \bar{u}=-\nabla \cdot[(u_{L}+\bar{u})\nabla(-\Delta)^{-1} (u_{L}+\bar{u})],\\
\bar{u}(x,0)=0.
\end{cases}
\end{equation}
Moreover, by using the Duhamel's principle, the system \eqref{eq2.2} can be further  reduced into the following equivalent integral equations:
\begin{equation}\label{eq2.3}
\bar{u}(t)=-\int_{0}^{t}e^{(t-s)\Delta}\nabla\cdot[(u_{L}+\bar{u})\nabla(-\Delta)^{-1} (u_{L}+\bar{u})]ds.
\end{equation}
Since $u_{L}$ is a solution of the heat equation with initial data $u_{0}$, one can easily obtain that for any $0<t\leq\infty$,
 \begin{equation}\label{eq2.1}
  \|u_{L}\|_{L^{1}_{t}(\dot{B}^{-2+\frac{d}{p}}_{p,1})}\leq \|u_{0}\|_{\dot{B}^{-2+\frac{d}{p}}_{p,1}}.
\end{equation}
Before giving the desired estimate of $\bar{u}$, we denote $f(t):= \|u_{L}(t)\|_{\dot{B}^{\frac{d}{p}}_{p,1}}$,
and for any positive real number $\lambda$, we introduce the weighted function
$\bar{u}_{\lambda,f}(x,t):=\bar{u}(x,t)\exp\{-\lambda\int_{0}^{t}f(\tau)d\tau\}$.
Back to \eqref{eq2.3},  it holds that
 \begin{align*}
  \bar{u}_{\lambda,f}&=-\int_{0}^{t}e^{(t-s)\Delta}e^{-\lambda\int_{0}^{t}f(\tau)d\tau}\nabla\cdot(u_{L}\nabla(-\Delta)^{-1} u_{L})ds\nonumber\\
  &-\int_{0}^{t}e^{(t-s)\Delta}e^{-\lambda\int_{s}^{t}f(\tau)d\tau}
  \nabla\cdot\big[u_{L}\nabla(-\Delta)^{-1} \bar{u}_{\lambda,f}+\bar{u}_{\lambda,f}\nabla(-\Delta)^{-1} u_{L}+\bar{u}\nabla (-\Delta)^{-1} \bar{u}_{\lambda,f}\big]ds.
 \end{align*}
 We know from the solvability of heat equation in the framework of the homogeneous Besov spaces that  there exists a positive constant $\kappa$ such that
 \begin{align}\label{eq2.5}
 \|&\bar{u}_{\lambda,f}\|_{\mathcal{L}^{\infty}_{t}(\dot{B}^{-2+\frac{d}{p}}_{p,1})}+\kappa\|\bar{u}_{\lambda,f}\|_{L^{1}_{t}(\dot{B}^{\frac{d}{p}}_{p,1})} \leq \|e^{-\lambda\int_{0}^{t}f(\tau)d\tau}u_{L}\nabla(-\Delta)^{-1} u_{L}\|_{L^{1}_{t}(\dot{B}^{-1+\frac{d}{p}}_{p,1})}\nonumber\\
 &+\|e^{-\lambda\int_{s}^{t}f(\tau)d\tau}\big[u_{L}\nabla(-\Delta)^{-1} \bar{u}_{\lambda,f}+\bar{u}_{\lambda,f}\nabla(-\Delta)^{-1} u_{L}+\bar{u}\nabla (-\Delta)^{-1} \bar{u}_{\lambda,f}\big]\|_{L^{1}_{t}(\dot{B}^{-1+\frac{d}{p}}_{p,1})}
 \nonumber\\
  &:=I_{1}+I_{2}+I_{3}+I_{4}.
\end{align}
 For $I_{1}$, one obtains that
 \begin{align}\label{eq2.6}
I_{1}\leq C \|u_{L}\nabla(-\Delta)^{-1} u_{L}\|_{L^{1}_{t}(\dot{B}^{-1+\frac{d}{p}}_{p,1})}.
 \end{align}
For $I_{2}$ and $I_{3}$, by using the following symmetric structure:
\begin{equation}\label{eq2.7}
    u\nabla(-\Delta)^{-1}v+v\nabla(-\Delta)^{-1}u=-\nabla\cdot\big(\nabla(-\Delta)^{-1}u\nabla(-\Delta)^{-1}v\big),
\end{equation}
we can deduce that
\begin{align}\label{eq2.8}
I_{2}+I_{3} &\leq C \int_{0}^{t}e^{-\lambda\int_{s}^{t}f(\tau)d\tau}\|\nabla\cdot(\nabla(-\Delta)^{-1}u_{L}\nabla(-\Delta)^{-1} \bar{u}_{\lambda,f})\|_{\dot{B}^{-1+\frac{d}{p}}_{p,1}}ds\nonumber\\
   &\leq C \int_{0}^{t}e^{-\lambda\int_{s}^{t}f(\tau)d\tau}\|\nabla(-\Delta)^{-1}u_{L}\nabla(-\Delta)^{-1} \bar{u}_{\lambda,f}\|_{\dot{B}^{\frac{d}{p}}_{p,1}}ds\nonumber\\
   &\leq C \int_{0}^{t}e^{-\lambda\int_{s}^{t}f(\tau)d\tau}\|\nabla(-\Delta)^{-1}u_{L}\|_{\dot{B}^{\frac{d}{p}}_{p,1}}\|\nabla(-\Delta)^{-1}\bar{u}_{\lambda,f}\|_{\dot{B}^{\frac{d}{p}}_{p,1}}ds\nonumber\\
     &\leq C \int_{0}^{t}e^{-\lambda\int_{s}^{t}f(\tau)d\tau}\|u_{L}\|_{\dot{B}^{-1+\frac{d}{p}}_{p,1}}\|\bar{u}_{\lambda,f}\|_{\dot{B}^{-1+\frac{d}{p}}_{p,1}}ds\nonumber\\
  &\leq C \int_{0}^{t}e^{-\lambda\int_{s}^{t}f(\tau)d\tau}\|u_{L}\|_{\dot{B}^{-2+\frac{d}{p}}_{p,1}}^{\frac{1}{2}}\|u_{L}\|_{\dot{B}^{\frac{d}{p}}_{p,1}}^{\frac{1}{2}}
  \|\bar{u}_{\lambda,f}\|_{\dot{B}^{-2+\frac{d}{p}}_{p,1}}^{\frac{1}{2}}\|\bar{u}_{\lambda,f}\|_{\dot{B}^{\frac{d}{p}}_{p,1}}^{\frac{1}{2}}ds\nonumber\\
    &\leq C\int_{0}^{t}e^{-\lambda\int_{s}^{t}f(\tau)d\tau}\big(\varepsilon \|\bar{u}_{\lambda,f}\|_{\dot{B}^{\frac{d}{p}}_{p,1}}+\|u_{L}\|_{\dot{B}^{-2+\frac{d}{p}}_{p,1}}\|u_{L}\|_{\dot{B}^{\frac{d}{p}}_{p,1}}
  \|\bar{u}_{\lambda,f}\|_{\dot{B}^{-2+\frac{d}{p}}_{p,1}}\big)ds\nonumber\\
   &\leq C \varepsilon \|\bar{u}_{\lambda,f}\|_{L^{1}_{t}(\dot{B}^{\frac{d}{p}}_{p,1})}+\|u_{0}\|_{\dot{B}^{-2+\frac{d}{p}}_{p,1}} \int_{0}^{t}e^{-\lambda\int_{s}^{t}f(\tau)d\tau} \|u_{L}\|_{\dot{B}^{\frac{d}{p}}_{p,1}}
  \|\bar{u}_{\lambda,f}\|_{\dot{B}^{-2+\frac{d}{p}}_{p,1}}ds,
\end{align}
where we have used the fact  that
 \begin{equation*}
  \|u_{L}\|_{L^{\infty}_{t}(\dot{B}^{-2+\frac{d}{p}}_{p,1})}\leq \|u_{L}\|_{\mathcal{L}^{\infty}_{t}(\dot{B}^{-2+\frac{d}{p}}_{p,1})}\leq \|u_{0}\|_{\dot{B}^{-2+\frac{d}{p}}_{p,1}}.
\end{equation*}
For $I_{4}$, since it has a nice structure as
$$
\bar{u}\nabla (-\Delta)^{-1} \bar{u}_{\lambda,f}=-\frac{1}{2}\nabla\cdot(\nabla (-\Delta)^{-1}\bar{u}\nabla (-\Delta)^{-1} \bar{u}_{\lambda,f}),
$$
we can estimate it as
\begin{align}\label{eq2.9}
 I_{4}
  &\leq C \int_{0}^{t}e^{-\lambda\int_{s}^{t}f(\tau)d\tau}\|\nabla (-\Delta)^{-1}\bar{u}\nabla (-\Delta)^{-1} \bar{u}_{\lambda,f}\|_{\dot{B}^{\frac{d}{p}}_{p,1}}ds\nonumber\\
  &\leq C \int_{0}^{t}\|\bar{u}\|_{\dot{B}^{-1+\frac{d}{p}}_{p,1}}\|\bar{u}_{\lambda,f}\|_{\dot{B}^{-1+\frac{d}{p}}_{p,1}}ds\nonumber\\
  &\leq C \int_{0}^{t}\|\bar{u}\|_{\dot{B}^{-2+\frac{d}{p}}_{p,1}}^{\frac{1}{2}}\|\bar{u}\|_{\dot{B}^{\frac{d}{p}}_{p,1}}^{\frac{1}{2}}
  \|\bar{u}_{\lambda,f}\|_{\dot{B}^{-2+\frac{d}{p}}_{p,1}}^{\frac{1}{2}}\|\bar{u}_{\lambda,f}\|_{\dot{B}^{\frac{d}{p}}_{p,1}}^{\frac{1}{2}}ds\nonumber\\
    &\leq C \int_{0}^{t}\|\bar{u}\|_{\dot{B}^{-2+\frac{d}{p}}_{p,1}} \|\bar{u}_{\lambda,f}\|_{\dot{B}^{\frac{d}{p}}_{p,1}}ds\nonumber\\
   &\leq C \|\bar{u}\|_{\mathcal{L}^{\infty}_{t}(\dot{B}^{-2+\frac{d}{p}}_{p,1})}\|\bar{u}_{\lambda,f}\|_{L^{1}_{t}(\dot{B}^{\frac{d}{p}}_{p,1})}.
\end{align}
Notice that
\begin{align*}
  \int_{0}^{t}e^{-\lambda\int_{s}^{t}f(\tau)d\tau}\|u_{L}(s)\|_{\dot{B}^{\frac{d}{p}}_{p,1}}ds=
  \int_{0}^{t}e^{-\lambda\int_{s}^{t}f(\tau)d\tau}f(s)ds
  =\frac{1}{\lambda}\int_{0}^{t}d\Big(e^{-\lambda\int_{s}^{t}f(\tau)d\tau}\Big)
  \leq \frac{1}{\lambda}.
\end{align*}
Thus taking all above estimates \eqref{eq2.6}--\eqref{eq2.9} into \eqref{eq2.5}, there exists a positive constant $C$ such that
\begin{align}\label{eq2.10}
\|&\bar{u}_{\lambda,f}\|_{\mathcal{L}^{\infty}_{t}(\dot{B}^{-2+\frac{d}{p}}_{p,1})} +\kappa\|\bar{u}_{\lambda,f}\|_{L^{1}_{t}(\dot{B}^{\frac{d}{p}}_{p,1})} \leq C\Big(\|u_{L}\nabla(-\Delta)^{-1} u_{L}\|_{L^{1}_{t}(\dot{B}^{-1+\frac{d}{p}}_{p,1})}\nonumber\\
 & +\frac{1}{\lambda}\|u_{0}\|_{\dot{B}^{-2+\frac{d}{p}}_{p,1}} \|\bar{u}_{\lambda,f}\|_{\mathcal{L}^{\infty}_{t}(\dot{B}^{-2+\frac{d}{p}}_{p,1})} +\varepsilon \|\bar{u}_{\lambda,f}\|_{L^{1}_{t}(\dot{B}^{\frac{d}{p}}_{p,1})} +\|\bar{u}\|_{\mathcal{L}^{\infty}_{t}(\dot{B}^{-2+\frac{d}{p}}_{p,1})}\|\bar{u}_{\lambda,f}\|_{L^{1}_{t}(\dot{B}^{\frac{d}{p}}_{p,1})}\Big).
\end{align}

Now we complete the proof of Theorem \ref{th1.1}. We know from \cite{I11} that there exists $T>0$ such that the system \eqref{eq1.1}  admits a unique local
solution $u\in \mathcal{L}^{\infty}(0,T;
  \dot{B}^{-2+\frac{d}{p}}_{p,1}(\mathbb{R}^{d}))\cap  L^{1}(0,T; \dot{B}^{\frac{d}{p}}_{p,1}(\mathbb{R}^{d}))$. Let us denote by $T_{*}$ the maximal existence time of this solution. Then to prove Theorem \ref{th1.1},
it suffices to prove $T_{*}=\infty$ under the initial condition \eqref{eq1.2}, and \eqref{eq1.4} holds. We prove it by contradiction. Assume that $T_{*}<\infty$, based on the estimate \eqref{eq2.10}, let $\eta$ be a small enough positive constant which the exact value will be determined later, we define $T_{\eta}$ as
\begin{align}\label{eq2.11}
  T_{\eta}:=\max\Big\{t\in[0,T_{*}):\ \|\bar{u}\|_{\mathcal{L}^{\infty}_{t}(\dot{B}^{-2+\frac{d}{p}}_{p,1})}
  +\kappa\|\bar{u}\|_{L^{1}_{t}(\dot{B}^{\frac{d}{p}}_{p,1})}\leq   \eta \Big\}.
\end{align}
Taking $\lambda$ large enough, $\varepsilon$ and $\eta$ small enough such that
$$
\lambda\geq2C\|u_{0}\|_{\dot{B}^{-2+\frac{d}{p}}_{p,1}},\ \ \ \varepsilon<\frac{\kappa}{4C}\ \ \text{and}\ \ \
\eta\leq \min\{1, \frac{\kappa}{4C}\},
$$
we can deduce from \eqref{eq2.10} to get that
\begin{align}\label{eq2.12}
   \|\bar{u}_{\lambda,f}\|_{\mathcal{L}^{\infty}_{t}(\dot{B}^{-2+\frac{d}{p}}_{p,1})}
   +\kappa\|\bar{u}_{\lambda,f}\|_{L^{1}_{t}(\dot{B}^{\frac{d}{p}}_{p,1})}
   &\leq 2C\|u_{L}\nabla(-\Delta)^{-1} u_{L}\|_{L^{1}_{t}(\dot{B}^{-1+\frac{d}{p}}_{p,1})}.
\end{align}
As a consequence,  for all $t\leq T_{\eta}$, it holds that
\begin{align}\label{eq2.13}
   \|\bar{u}\|_{\mathcal{L}^{\infty}_{t}(\dot{B}^{-2+\frac{d}{p}}_{p,1})}
   +\kappa\|\bar{u}\|_{L^{1}_{t}(\dot{B}^{\frac{d}{p}}_{p,1})}
      \leq 2C\|u_{L}\nabla(-\Delta)^{-1} u_{L}\|_{L^{1}_{t}(\dot{B}^{-1+\frac{d}{p}}_{p,1})}
  \exp\Big\{\lambda\int_{0}^{t}\|u_{L}(\tau)\|_{\dot{B}^{\frac{d}{p}}_{p,1}}d\tau\Big\}.
\end{align}
Based on \eqref{eq2.1}, we obtain that there exists a positive constant $C$ such that for all $t\leq T_{\eta}$, it holds that
\begin{align}\label{eq2.15}
   \|\bar{u}\|_{\mathcal{L}^{\infty}_{t}(\dot{B}^{-2+\frac{d}{p}}_{p,1})}
   +\kappa\|\bar{u}\|_{L^{1}_{t}(\dot{B}^{\frac{d}{p}}_{p,1})}
      \leq C\|u_{L}\nabla(-\Delta)^{-1} u_{L}\|_{L^{1}_{t}(\dot{B}^{-1+\frac{d}{p}}_{p,1})}
   \exp\Big\{C\|u_{0}\|_{\dot{B}^{-2+\frac{d}{p}}_{p,1}}\Big\}.
\end{align}
We finally conclude that if we take $C_{0}$ large enough and $c_{0}$ small enough in \eqref{eq1.2}, then it follows from \eqref{eq2.15} that
\begin{align*}
  \|\bar{u}\|_{\mathcal{L}^{\infty}_{t}(\dot{B}^{-2+\frac{d}{p}}_{p,1})}
   +\kappa\|\bar{u}\|_{L^{1}_{t}(\dot{B}^{\frac{d}{p}}_{p,1})}
  \leq
  \frac{\eta}{2}
\end{align*}
for all $t<T_{\eta}$, which contradicts with the maximality of $T_{\eta}$, thus $T^*=\infty$, and \eqref{eq1.4} holds true. We complete the proof of Theorem \ref{th1.1}.

\section{Proof of Theorem \ref{th1.2}}

In this section we give the proof of Theorem \ref{th1.2}. We shall check that the initial data $u_{0,\varepsilon}$ introduced in Theorem \ref{th1.2} satisfies the nonlinear smallness assumption \eqref{eq1.2}, but its $\dot{B}^{-2}_{\infty,\infty}$ norm is greater than $(-\log \varepsilon)^{\frac{1}{5}}$.
\begin{lemma}\label{le3.1}{\em (see \cite{CG09})}
Let $f\in\mathcal{S}(\mathbb{R}^3)$ be given and $\sigma\in(0, 3-\frac{3}{p})$. Then there
exists a constant $C$ such that for any $\varepsilon, \alpha\in(0,1)$,  the function $f_{\varepsilon}(x):=e^{i\frac{x_{3}}{\varepsilon}}f(x_{1},\frac{x_{2}}{\varepsilon^{\alpha}},x_{3})$ satisfies,
for all $1\leq p\leq \infty$,
\begin{equation}\label{eq3.1}
 \|f_{\varepsilon}\|_{\dot{B}^{-\sigma}_{p,1}}\leq C\varepsilon^{\sigma+\frac{\alpha}{p}} \ \ \text{and}\ \ \|f_{\varepsilon}\|_{\dot{B}^{-\sigma}_{\infty,\infty}}\geq C^{-1}\varepsilon^{\sigma}.
\end{equation}
\end{lemma}

\begin{lemma}\label{le3.2} {\em(see \cite{L18})}
Let $f$ and $g$ be in $\dot{B}^{-1}_{2p,1}(\mathbb{R}^3)\cap \dot{H}^{-1}(\mathbb{R}^3)$. Then for all $3<p\leq \infty$,  we have
\begin{equation}\label{eq3.2}
   \| e^{t\Delta}f e^{t\Delta}g\|_{L^{1}(\mathbb{R}^{+}; \dot{B}^{-1+\frac{3}{p}}_{p,1})}\leq C\big(\|f\|_{\dot{B}^{-1}_{2p,1}}\|g\|_{\dot{B}^{-1}_{2p,1}}\big)^{\frac{2p}{3p-3}}
  \big(\|f\|_{\dot{H}^{-1}}\|g\|_{\dot{H}^{-1}}\big)^{\frac{p-3}{3p-3}}.
\end{equation}
\end{lemma}

Define a smooth function as $u_{0,\varepsilon}:=\nabla\cdot w_{0,\varepsilon}$, where  $w_{0,\varepsilon}$ is defined as
$w_{0,\varepsilon}:=(\partial_{2}\varphi_{\varepsilon}, \partial_{1}\varphi_{\varepsilon}, 0)$ and
 \begin{equation*}
 \varphi_{\varepsilon}(x):=\frac{(-\log \varepsilon)^{\frac{1}{5}}}{\varepsilon^{1-\alpha}}\cos(\frac{x_{3}}{\varepsilon})\phi(x_{1}, \frac{x_{2}}{\varepsilon^{\alpha}},x_{3}).
\end{equation*}
Then we have the following upper and lower bound of $w_{0,\varepsilon}$ in Besov type space.
\begin{lemma}\label{le3.3} {\em(see \cite{CG09})}
A constant $C$ exists such that, for all $p\geq 1$, we have
\begin{equation}\label{eq3.3}
   \|w_{0, \varepsilon}\|_{\dot{B}^{-1}_{p,1}}\leq  C\varepsilon^{\frac{\alpha}{p}}(-\log \varepsilon)^{\frac{1}{5}}\ \ \text{and}\ \ \|w_{0, \varepsilon}\|_{\dot{B}^{-1}_{\infty,\infty}}\geq  C^{-1}(-\log \varepsilon)^{\frac{1}{5}}.
\end{equation}
\end{lemma}

By using Lemma \ref{le3.1}, we know that there exists a positive constant $C$ such that
\begin{equation*}
\|u_{0, \varepsilon}\|_{\dot{B}^{-2}_{\infty,\infty}}=\|\nabla\cdot w_{0, \varepsilon}\|_{\dot{B}^{-2}_{\infty,\infty}}\approx\|w_{0, \varepsilon}\|_{\dot{B}^{-1}_{\infty,\infty}}\geq  C^{-1}(-\log \varepsilon)^{\frac{1}{5}}.
\end{equation*}
Then for all $1\leq p<\infty$,  it follows the Sobolev embedding result that
\begin{equation}\label{eq3.3}
  \|u_{0, \varepsilon}\|_{\dot{B}^{-2+\frac{3}{p}}_{p,1}}\geq \|u_{0, \varepsilon}\|_{\dot{B}^{-2}_{\infty,\infty}}\geq C^{-1}(-\log \varepsilon)^{\frac{1}{5}}.
\end{equation}
On the other hand, it can easily calculate
%
%
%
\begin{align*}
  \nabla\cdot w_{0,\varepsilon}=\frac{2}{\varepsilon}(-\log \varepsilon)^{\frac{1}{5}}\cos(\frac{x_{3}}{\varepsilon})\partial_{1}\partial_{2}\phi(x_{1}, \frac{x_{2}}{\varepsilon^{\alpha}},x_{3}),
\end{align*}
which yields that
\begin{equation}\label{eq3.5}
  e^{t\Delta}w_{0, \varepsilon}\nabla \cdot e^{t\Delta}w_{0, \varepsilon}=\frac{2}{\varepsilon^{2-\alpha}}(-\log \varepsilon)^{\frac{2}{5}}e^{t\Delta} f_{\varepsilon}e^{t\Delta} g_{\varepsilon}
  +\frac{2}{\varepsilon^{2}}(-\log \varepsilon)^{\frac{2}{5}}e^{t\Delta} \tilde{f}_{\varepsilon}e^{t\Delta} \tilde{g}_{\varepsilon},
\end{equation}
where $ f_{\varepsilon}, g_{\varepsilon}$, $\tilde{ f}_{\varepsilon}$ and $\tilde{g}_{\varepsilon}$ are smooth functions. By using Lemmas \ref{le3.1} and \ref{le3.2}, it follows that for all $3<p\leq \infty$, we have
\begin{align}\label{eq3.6}
  \| e^{t\Delta} f_{\varepsilon}e^{t\Delta} g_{\varepsilon}\|_{L^{1}(\mathbb{R}^{+};\dot{B}^{-1+\frac{3}{p}}_{p,1})}&\leq C\big(\|f_{\varepsilon}\|_{\dot{B}^{-1}_{2p,1}}\|g_{\varepsilon}\|_{\dot{B}^{-1}_{2p,1}}\big)^{\frac{2p}{3p-3}}
  \big(\|f_{\varepsilon}\|_{\dot{H}^{-1}}\|g_{\varepsilon}\|_{\dot{H}^{-1}}\big)^{\frac{p-3}{3p-3}}\nonumber\\
  &\leq C\big(\|f_{\varepsilon}\|_{\dot{B}^{-1}_{2p,1}}\|g_{\varepsilon}\|_{\dot{B}^{-1}_{2p,1}}\big)^{\frac{2p}{3p-3}}
  \big(\|f_{\varepsilon}\|_{\dot{B}^{-1}_{2,1}}\|g_{\varepsilon}\|_{\dot{B}^{-1}_{2,1}}\big)^{\frac{p-3}{3p-3}}\nonumber\\
  &\leq C\varepsilon^{2+\frac{\alpha}{3}}.
\end{align}
Similarly,
\begin{align}\label{eq3.61}
  \| e^{t\Delta} \tilde{f}_{\varepsilon}e^{t\Delta} \tilde{g}_{\varepsilon}\|_{L^{1}(\mathbb{R}^{+};\dot{B}^{-1+\frac{3}{p}}_{p,1})}\leq C\varepsilon^{2+\frac{\alpha}{3}}.
\end{align}
Combining \eqref{eq3.5},  \eqref{eq3.6} and \eqref{eq3.61}, we immediately obtain
\begin{equation*}
  \|e^{t\Delta}w_{0, \varepsilon}\nabla \cdot e^{t\Delta}w_{0, \varepsilon}\|_{L^{1}(\mathbb{R}^{+};\dot{B}^{-1+\frac{3}{p}}_{p,1})}\leq C\varepsilon^{\frac{\alpha}{3}}(-\log \varepsilon)^{\frac{2}{5}}.
\end{equation*}
Notice that $u_{0,\varepsilon}=\nabla\cdot w_{0,\varepsilon}$,  which gives us to $w_{0,\varepsilon}=-\nabla(-\Delta)^{-1} u_{0,\varepsilon}$, thus we finally obtain
\begin{equation*}
  \|e^{t\Delta}u_{0, \varepsilon}\nabla(-\Delta)^{-1}e^{t\Delta}u_{0, \varepsilon}\|_{L^{1}(\mathbb{R}^{+};\dot{B}^{-1+\frac{3}{p}}_{p,1})}\leq C\varepsilon^{\frac{\alpha}{3}}(-\log \varepsilon)^{\frac{2}{5}}.
\end{equation*}
We complete the proof of Theorem \ref{th1.2}.

\small{

}

\end{document}